\input amstex
\mag=\magstep1
\documentstyle{amsppt}
\nologo
\NoBlackBoxes
% \NoPageNumbers
\NoRunningHeads
\pagewidth{360pt}
\pageheight{580pt}
%\pageheight{600pt} %  for arXiv.org
%\vcorrection{-25pt} % for arXiv.org
% \leftskip=-7mm
\define\dint{\dsize\int}
\def\fp{\flushpar}
\font\tenptmbit=cmmib10 % Computer Modern Mathematical italic bold of 10pt
\font\eightptmbit=cmmib8
\define\bk#1{\text{\tenptmbit #1}}
\define\eightptbk#1{\text{\eightptmbit{#1}}}
\define\inbox#1{$\boxed{\text{#1}}$}
\define\underbarl#1{\lower 1.4pt \hbox{\underbar{\raise 1.4pt \hbox{#1}}}}
\define\ce#1{\lceil#1\rceil} 
\define\tp#1{\negthinspace\ ^t#1}
\define\dg#1{(d^{\circ}\geqq#1)}
\define\Dg#1#2{(d^{\circ}(#1)\geqq#2)}

\define\jrac#1#2{\dfrac{\lower 2pt\hbox{$#1$}}{\raise 2pt\hbox{$#2$}}}
\define\sz#1{^{\raise 1pt\hbox{$\sssize #1$}}\hskip -1pt}
\define\lr#1{^{\sssize\left(#1\right)}}
\define\br#1{^{\sssize\left[#1\right]}}
\define\do#1{^{\sssize\left<#1\right>}}

\define\rank{\text{rank}}
\define\twolower#1{\hbox{\raise -2pt \hbox{$#1$}}}

\define\hookdownarrow{\kern 0pt
 \hbox{\vbox{\offinterlineskip \kern 0pt
  \hbox{$\cap$}\kern 0pt
  \hbox{\hskip 3.5pt$\downarrow$}\kern 0pt}}
}  % $B#d(Bvi $B$G$O%:%l%k!$(BPS $B$J$iB7$U(B. 

\def\maru#1{{\ooalign{\hfil#1\/\hfil\crcr\raise.167ex\hbox{\mathhexbox20D}}}} 
\font\tenptmbit=cmmib10 % Computer Modern Mathematical italic bold of 10pt
\font\eightptmbit=cmmib8

% Type-writer type fonts

% Computer Modern Bold Face

\font\thirteenptbf=cmbx8 scaled \magstep3
\font\twelveptbf=cmbx12

% Computer Modern Mathematical Italic
 
% Computer Modern Roman

\font\twelveptcm=cmr12 

\font\sixptcm=cmr6
% Computer Modern ??? Sub-Script 

%
\redefine\jtag#1{\tag"$B!J(B#1$B!K(B"}
\def\boxit#1{\vbox{\hrule\hbox{\vrule\kern3pt
    \vbox to 43pt{\hsize 182pt\kern3pt#1\eject\kern3pt\vfill}
    \kern1pt\vrule}\hrule}} % $B$3$l$,4pK\7?(B
\def\oVrule{\vrule width .5pt}
\def\oHrule{\hrule height .5pt}
%%%%%%%%%% $B%Q%i%0%i%UMQ$N%\%C%/%9(B %%%%%%%%%%%%%%%%%%%%%%
%%%%%%%%%% $B2#(B(pt)$B!$=D(B(pt)$B!$FbMF(B %%%%%%%%%%%%%%%%%%%%%%%%
%%% $B7S@~$N6R$,(B .5pt$B!$(B $B>e2<M>Gr$,(B 8pt $B:81&$NM>Gr$,(B 10pt $B$E$D$N(B
%%% $B%\%C%/%9$,$G$-$k(B. $B=D2#$ND9$5$O!$7S@~$^$G4^$a$?CM$rF~NO$;$h(B. 
\def\lbox#1#2#3{\kern 0pt % $B$3$N(B \kern $B@_Dj$,$J$$$H:8$K7d4V$,$G$-$k(B. 
\dimen1=#1pt \dimen2=#2pt% 
\advance\dimen1 by -1.0pt% $B7S@~$NB@$5$N#2G\$r:9$70z$/(B.  
\dimen3=\dimen1% 
\advance\dimen3 by -16pt% $BFbMF$N>eC<$+$i7S@~2<C<$^$G$N5wN%(B(=$BM>Gr(B)$B$N#2G\$r:9$70z$/(B
\advance\dimen2 by -1.0pt% $B7S@~$NB@$5$N#2G\$r:9$70z$/(B.  
\dimen4=\dimen2% 
\advance\dimen4 by -20pt% $BFbMF$N:8C<$+$i7S@~1&C<$^$G$N5wN%(B(=$BM>Gr(B)$B$N#2G\$r:9$70z$/(B
\vbox to #1pt%      <---------------------- #2 --------------------->             
{%                  +-----------------------------------------------+             
\hsize #2pt%        |########### height 0.5pt ######################|             
\oHrule%            |---+--------------------------------+------+---+  -----------
\hbox to #2pt%      |###|          Here \vskip 8pt (A)          |###|         ^   
{%                  |# #+------+-------------------------+------+# #|  -----  |   
\vsize \dimen1%     |#w#|      |/////////////////////////|      |#w#|    ^    |   
\oVrule%            |#i#|      |/////////////////////////|      |#i#|    |    |   
\vbox to \dimen1%   |#d#|Here  |//HERE/IS/THE/CONTENTS///|Here  |#d#|    |    |   
{%                  |#t#|\hskip|////////(#3)/////////////|\hskip|#t#|    | \dimen1
\hsize \dimen2%     |#h#|10pt  |/////////////////////////|10pt  |#h#|    |    |   
\vskip 8pt%   (A)   |# #| (C)  |<----- \dimen4 --------->| (D)  |# #| \dimen3 |   
\hbox to \dimen2%   |#0#|      |/////////////////////////|      |#0#|    |    |   
{%                  |#.#|      |/////////////////////////|      |#.#|    |    |   
%                   |#5#|      |/////////////////////////|      |#5#|    v    |   
%                   |#p#|      |/////////////////////////|      |#p#|         |   
%                   |#t#+------+-------------------------+------+#t#|  -----  |   
\vsize \dimen3%     |###|          Here \vskip 8pt (B)          |###|         v   
\hskip 10pt% ...(C) +---+---------------------------------------+---+  ---------- 
%                   |######### height 0.5 pt #######################|             
%                   +-----------------------------------------------+             
%                                                                                 
%                       |<------------\dimen2 ----------------->|                 
%                                                                                 
\vbox to \dimen3{\hsize \dimen4\fp #3\vfil}%
% \vbox{\hrule width\dimen4 height\dimen3}% $B%F%9%HMQ$N$L$j$D$V$7(B. 
\hfil\hskip 10pt% ...(D) 
}%
\vskip 8pt%     ....(B)
}%
\oVrule%
}%
\oHrule%
}\kern 0pt% $B$3$l$,$J$$$H7d4V$,3+$$$F$7$^$U(B. 
}%
%%%%%%%%%% $BFbIt$O%;%s%?%j%s%0$9$k%\%C%/%9(B %%%%%%%%%%%%%%%%%%%%%
%%%%%%%%%% $B=D(B(pt)$B!$2#(B(pt)$B!$FbMF(B %%%%%%%%%%%%%%%%%%%%%%%%
%%% $B7S@~$N6R$,(B .5pt$B!$(B $B>e2<:81&$NM>Gr$,$=$l$>$l(B 1.5pt $B$E$D$N(B
%%% $B%\%C%/%9$,$G$-$k(B. $B=D2#$ND9$5$O!$7S@~$^$G4^$a$?CM$rF~NO$;$h(B. 
\def\cbox#1#2#3{\kern 0pt % $B$3$N(B \kern $B@_Dj$,$J$$$H:8$K7d4V$,$G$-$k(B. 
\dimen1=#1pt \dimen2=#2pt% 
\advance\dimen1 by -1.0pt% $B7S@~$NB@$5$N#2G\$r:9$70z$/(B.  
\dimen3=\dimen1% 
\advance\dimen3 by -3pt% $BFbMF$N>eC<$+$i7S@~2<C<$^$G$N5wN%(B(=$BM>Gr(B)$B$N#2G\$r:9$70z$/(B
\advance\dimen2 by -1.0pt% $B7S@~$NB@$5$N#2G\$r:9$70z$/(B.  
\dimen4=\dimen2% 
\advance\dimen4 by -3pt% $BFbMF$N:8C<$+$i7S@~1&C<$^$G$N5wN%(B(=$BM>Gr(B)$B$N#2G\$r:9$70z$/(B
\vbox to #1pt%
{%
\hsize #2pt
\oHrule%
\hbox to #2pt%
{%
\vsize \dimen1%
\oVrule%
\vbox to \dimen1%
{%
\hsize \dimen2%
\vskip 1.5pt%  ...(A)
\hbox to \dimen2%
{%
\vsize \dimen3%
\hskip 1.5pt\hfil%  ...(C)
\vbox to \dimen3{\hsize \dimen4\vfil\hbox{#3}\vfil}%
% \vbox{\hrule width\dimen4 height\dimen3}% $B%F%9%HMQ$N$L$j$D$V$7(B. 
\hfil\hskip 1.5pt%  ...(D)
}%
\vskip 1.5pt%       ...(B)
}%
\oVrule%
}%
\oHrule%
}\kern 0pt% $B$3$l$,$J$$$H7d4V$,3+$$$F$7$^$U(B. 
}%
\def\cboxit#1#2#3{$\hbox{\lower 2.5pt \hbox{\cbox{#1}{#2}{#3}}}$}

\document
\baselineskip=13.3pt % 2\baselineskip % 24pt
{\ }
\vskip 5pt
\centerline{\thirteenptbf Abelian Functions for Trigonal Curves of}
\vskip 3pt
\centerline{\thirteenptbf Degree Four and Determinantal Formulae}
\vskip 3pt
\centerline{\thirteenptbf in Purely Trigonal Case}
\vskip 12pt
\centerline{\twelveptcm Yoshihiro \^ONISHI}
\vskip 39pt
\fp

\fp
{\twelveptbf Preface}
\vskip 3pt
\fp
In the theory of elliptic functions, 
there are two kinds of formulae of Frobenius-Stickelberger \cite{{\bf FS}} and of Kiepert \cite{{\bf K}}, 
both of which connect the function  $\sigma(u)$  with  $\wp(u)$  and its (higher) derivatives 
through an determinantal expression.  
These formulae were naturally generalized to hyperelliptic functions 
by the papers \cite{{\bf \^O2}}, \cite{{\bf \^O3}}, and \cite{{\bf \^O4}}.  
Avoiding an unnecessary generality, we restrict the story only for the 
simplest purely trigonal curve  $y^3=x^4+\cdots$,
where the right hand side is a monic biquadratic polynomial of  $x$  
(Theorem 4.3 and Corollary 5.2).  
For more general purely trigonal curve, 
or for any purely $d$-gonal curve ($d=4$, $5$, $\cdots$), 
the author would like to publish in other papers, 
including formulae of Cantor-type (see \cite{{\bf \^O4}}).  

In the case of hyperelliptic functions, 
we considered only the hyperellptic curves ramified at infinity.  
The theta divisor of the Jacobian variety of such a curve is symmetric 
with respect to the origin of the Jacobian variety.   
Each of purely trigonal curves considered in this paper 
is also completely ramified at infinity 
and it is acted by the third roots of unity.  
Hence, the theta divisor of the Jacobian variety of a purely trigonal curve  
has third order symmetry with respect to the origin.  
Similar symmetry is also possessed by any purely $d$-gonal curve.  

On the other hand, any curve that is not purely $d$-gonal 
does not have symmetry with respect to the origin at all. 
Since this fact is very serious, the author never imagine 
whether generalizations of Frobenius-Stickelberger type and 
of Kiepert type exist or not.  

To work out this paper was very simplified by virtue of the three nice papers 
\cite{{\bf BEL1}}, \cite{{\bf BEL2}}.  
The author would like to express hearty thanks to the three of authors of 
these papers.  
\vskip 7pt

\newpage
\fp
{\bf Conventions}.  Needless to say, but we denote 
the ring of integers by  $\bold Z$, 
the field of real numbers by  $\bold R$, 
the field of complex numbers by  $\bold C$.  
The transpose of a vector  $u$  is denoted by  $\tp{u}$.  
The symbol  $d^{\circ}(z_1,\cdots,z_m)\geqq d$  means 
a power series whose all terms with respect to 
the specified variables  $z_1$, $\cdots$, $z_m$  are 
of total degree bigger than  $d$.  
Please note that this symbol does not means 
it is a power series which contains 
only the variables  $z_1$, $\cdots$, $z_m$.  
\vskip 30pt
\fp
\hskip 10pt{\twelveptbf Contents} 
\vskip 4pt
\fp\hskip 10pt
1. Preliminaries\fp\hskip 10pt
2. The sigma function (General case) \fp\hskip 10pt
3. The sigma function (Purely trigonal case) \fp\hskip 10pt
4. Frobenius-Stickelberbger-type formulae \fp\hskip 10pt
5. Kiepert-type formula 

\newpage
\fp
{\twelveptbf 1. Preliminaries}
\vskip 3pt
\fp
This and the next Sections are devoted for preliminaries under 
more general situation than the Sections later.   
Let  $C$  be a complete projective algebraic curve defined by  
  $$
  \aligned
  f(x,y)&=0, \ \ \text{where} \\
  f(x,&y)= y^3 - (\lambda_1 x + \lambda_4)y^2 - (\lambda_2 x^2 + \lambda_5 x + \lambda_8)y \\
  &-(x^4 + \lambda_3x^3 + \lambda_6x^2 + \lambda_9x + \lambda_{12}) \ \ \text{($\lambda_j$  are constants)}, 
  \endaligned
  \tag{1.1}
  $$
with the unique point   $\infty$  at infinity.  
The genus of  $C$  is  3  if  $C$  is non-singular.  
The set of three forms of the first kind (namely, holomorphic forms) 
  $$
  \omega_1=\frac{dx}{\frac{\partial}{\partial y}f(x,y)},  \ \ 
  \omega_2=\frac{xdx}{\frac{\partial}{\partial y}f(x,y)},  \ \  
  \omega_3=\frac{ydx}{\frac{\partial}{\partial y}f(x,y)}
  \tag{1.2}
  $$
forms a basis of the space of holomorphic 1-forms.  
We denote the vector whose coordinates are (1.2) by
  $$
  \omega=(\omega_1, \omega_2, \omega_3).
  \tag{1.3}
  $$
The general theory of Abelian integrals shows that the integrals 
  $$
  \aligned
  u&=(u_1,u_2,u_3) \\
  &=\int_{\infty}^{(x_1,y_1)}\omega
   +\int_{\infty}^{(x_2,y_2)}\omega
   +\int_{\infty}^{(x_3,y_3)}\omega
  \endaligned
  \tag{1.4}
  $$
with respect to all the path from  $\infty$  to 
three variable points  $(x_1,y_1)$, $(x_2,y_2)$, $(x_3,y_3)$  on  $C$  fill 
the whole space  $\bold C^3$.  
In this paper, we denote by the letters  $u$, $v$, $u\lr{j}$  various points in  $\bold C^3$,  
and the same letters with subscripts  $(u_1,u_2,u_3)$, $(v_1,v_2,v_3)$, $(u\lr{j}_1,u\lr{j}_2,u\lr{j}_3)$  
denote their canonical coordinates of  $\bold C^3$.  
We denote by  $\Lambda$  all the values of the integral above with respect to all the closed paths.    
Then  $\Lambda$  is a lattice of  $\bold C^3$.  
The  $\bold C$-valued points of the Jacobian variety of  $C$  is  $\bold C^3/\Lambda$.  
We denote this by  $J$.  
The canonical map given by modulo  $\Lambda$  is denoted by  $\kappa$:
  $$
  \kappa:\bold C^3\rightarrow \bold C^3/\Lambda=J.
  \tag{1.5}
  $$
Obviously, $\Lambda=\kappa^{-1}\big((0,0,0)\big)$.  
We have the standard embedding of  $C$  into  $J$  given by
  $$
  \aligned
  \iota :\ &C\hookrightarrow J \\
  &P \mapsto \int_{\infty}^P\hskip -3pt\omega \ \ \hbox{mod}\ \Lambda.
  \endaligned
  \tag{1.6}
  $$
Then  $\kappa^{-1}\iota(C)$  is a {\it universal Abelian covering} of  $C$.  
More generally, for  $0\leqq k\leqq 3$,  
the image of the  $k^{\text{\sixptcm th}}$  symmetric product  $\text{Sym}^k(C)$  of   $C$
by the map
  $$
  \aligned
  \iota :\ \text{Sym}^k(C)&\rightarrow J \\
  (P_1,\cdots,P_k) &\mapsto 
  \Big(\int_{\infty}^{P_1}\hskip -5pt\omega+\cdots+\int_{\infty}^{P_k}\hskip -5pt\omega\Big) \ 
  \hbox{\rm mod}\ \Lambda
  \endaligned
  \tag{1.7}
  $$
is denoted by  $W\br{k}$.  
Especially, $W\br{0}=(0,0,0)$, $W\br{1}=\iota(C)$, $W\br{3}=J$.  
We denote by  $\lceil -1\rceil$  the multiplication by  $(-1)$  in  $J$, namely
  $$
  {\lceil -1\rceil}(u_1,u_2,u_3)=(-u_1,-u_2,-u_3).
  \tag{1.8}
  $$
We define by this operation  ${\lceil -1\rceil}$  that
  $$
  \Theta\br{k}=W\br{k}\cup {\lceil -1\rceil}W\br{k}.
  \tag{1.9} 
  $$
We call  $\Theta\br{k}$  the {\it standard theta subvarieties}.  
Especially, we have   $\Theta\br{0}=(0,0,0)$  and  $\Theta\br{3}=J$.  
Note that we have  $\Theta\br{1}\neq W\br{1}$, $\Theta\br{2}\neq W\br{2}$  
as our case is different from the case of hyperelliptic curves.  
\vskip 5pt
\fp
\lbox{62}{360}
% --------------------------------------------------------------------------------------s
{{\bf  Lemma 1.10} \it \ Suppose  $u=(u_1,u_2,u_3)\in\kappa^{-1}\iota(C)$.  
Then  $u_1$  and  $u_2$  are expanded as a power series with respect to  $u_3$  of the forms
  $$
  u_1=\tfrac15{u_3}^5+\cdots,\ \ \ \ 
  u_2=\tfrac12{u_3}^2+\cdots.
  $$
}%---------------------------------------------------------------------------------------------e
\vskip 5pt
\fp
{\it Proof}. \ 
Taking a local parameter  $t=1/\root 3\of{x}$  on  $C$  at  $\infty$, 
we compute the integral
  $$
  u_j=\int_{\infty}^{(x,y)}du_j,
  \tag{1.11}
  $$
then we have 
  $$
  u_1=\tfrac15t^5+\cdots,\ \ \ \ 
  u_2=\tfrac12t^2+\cdots,\ \ \ \ 
  u_3=t+\cdots.
  \tag{1.12}
  $$
Hence the statement. 
\qed  
\vskip 10pt
\fp
The result above shows that  $u_3$  is a {\it local parameter} on  $\kappa^{-1}\iota(C)$  at the point  $(0,0,0)$.  
Now we consider the integral
  $$
  u=(u_1,u_2,u_3)=\int_{\infty}^{(x,y)}\omega,
  \tag{1.13}
  $$
and we denote its inverse function by 
  $$
  u\mapsto \big(x(u),y(u)\big).
  \tag{1.14}
  $$
\vskip 5pt
\fp
\lbox{69}{360}%--------------------------------------------------------------------------s
{{\bf Lemma 1.15}\it \  If  $u=(u_1,u_2,u_3)\in\kappa^{-1}\iota(C)$,   
then  $x(u)$  and  $y(u)$  are expanded as power series with respect to  $u_3$  as follows{\rm :}
  $$
  x(u)=\frac1{{u_3}^3}+\cdots,\ \ \ \ 
  y(u)=\frac1{{u_3}^4}+\cdots.
  $$
}%---------------------------------------------------------------------------------------------e
\vskip 5pt
\fp
{\it Proof}. \ 
This is proved similarly with Lemma 1.10.  
\qed
\vskip 10pt
\fp
\lbox{65}{360}{{\bf Definition 1.16} \it 
We define a weight called {\rm Sato weight} for appeared constants and variables 
as follows. 
The Sato weight of variables  $u_1$, $u_2$, $u_3$  are $5$, $2$, $1$, respectively, 
the Sato weight of each the coefficient  $\lambda_j$  in {\rm (1.1)} is  $-j$, 
The Sato weight of  $x(u)$  and  $y(u)$  are  $-4$  and  $-3$, respectively.  
}
\vskip 3pt
\fp
Under this convention, the formulae in this paper are of homogeneous weight. 
\def\rslt{\hbox{\rm rslt}}

We define the {\it discriminant} of  $C$  by
  $$
  D=\big[\rslt_x
  \big(\rslt_y\big(f(x,y), \tfrac{\partial}{\partial y}f(x,y)\big), 
       \rslt_y\big(f(x,y), \tfrac{\partial}{\partial y}f(x,y)\big)
  \big)\big]^{1/3},
  \tag{1.17}
  $$
where  $\rslt_z$  is the resultant of Sylvester with respect to the variable  $z$, 
and the cubic root is taken such as 
  $$
  D\in \bold Z[\lambda_1,\lambda_2,\lambda_3,\lambda_4,
    \lambda_5,\lambda_6,\lambda_8,\lambda_9,\lambda_{12}]
  \tag{1.18}
  $$
by regarding  $\lambda_j$s are indeterminates.   
While this possibility is not obvious, we omit it because of less importance in this paper. 

\newpage  

\fp
{\twelveptbf 2. Sigma function (General case)} 
\vskip 5pt
\fp
We now define the {\it sigma function}
  $$
  \sigma(u)=\sigma(u_1,u_2,u_3)
  \tag{2.1}
  $$
associating to  $C$  in the following \cite{{\bf Ba}}.  
We choose a set of generators 
  $$
  \alpha_i, \ \ \alpha_j \ (1\leqq i,\ j\leqq 3)
  \tag{2.2}
  $$
of  $H_1(C,\bold Z)$  such that their intersection numbers are 
$\alpha_i\cdot\alpha_j=\beta_i\cdot\beta_j=\delta_{ij}$  and  $\alpha_i\cdot\beta_j=0$.  
We denotes the period matrix obtained from the differentials (1.2) by
  $$
  [\,\omega'\ \ \omega'']=\bigg[\int_{\alpha_i}\omega_j \ \ \int_{\beta_i}\omega_j\bigg]_{i,j=1,2,3}.
  \tag{2.3}
  $$
For indeterminates  $X$, $Y$, $Z$, and $W$, we consider that 
  $$
  \Omega\big((X,Y),(Z,W)\big)=\frac{1}{(X-Z)\frac{\partial}{\partial Y}f(X,Y)}
  \sum_{k=1}^3Y^{3-k}\bigg[\frac{f(Z,W)}{W^{k-1}}\bigg]_W, 
  \tag{2.4}
  $$
where $[\ \ ]_W$  means taking only the terms of non-negative powers with respect to  $w$.  

\vskip 3pt
\fp
\lbox{280}{360}
{{\bf Lemma-Definition 2.5} \ {\rm (fundamental 2-forms of second kind)} \ 
\it Let
  $$
  \big((x,y),(z,w)\big)\mapsto R\big((z,w),(x,y)\big)dzdx
  \tag{2.6}
  $$
be a \,$2$-form on  $C\times C$  satisfying
  $$
  \lim_{x\to z}(z-x)^2R\big((z,w),(x,y)\big)=1, 
  \tag{2.7}
  $$
having only poles along the points such that  $(x,y)$\,$=$\,$(z,w)$,  
and being holomorphic elsewhere. 
For the differential {\rm (1.2)} of first kind, for  $\Omega$  in {\rm (2.4)}, 
and for two variables  $(x,y)$, $(z,w)$  on  $C$,  
there exist differential forms   $\eta_j=\eta_j(x,y)$ {\rm (}$j=1$, $2$, $3${\rm )} of second kind 
that have poles only at  $\infty$, for which such a $2$-form as above is expressed in the form
  $$
  R\big((x,y),(z,w)\big)
  :=\frac{d}{dx}\Omega\big((x,y),(z,w)\big)
  +\sum_{j=1}^3\frac{\omega_j(x, y)}{dx}\frac{\eta_j(z, w)}{dz},
  \tag{2.8}
  $$
where the derivation\footnotemark is one with respect to  $(x,y)\in C$.
We further require that  $R\big((x,y),(z,w)\big)$  is 
of homogeneous Sato weight {\rm(}hence weight $6${\rm)}, 
and the symmetricity 
  $$
  R\big((z,w),(x,y)\big)=R\big((x,y),(z,w)\big).
  \tag{2.9}
  $$
Such the set  $\{\eta_j\}$  exists and is uniquely determined 
modulo the space spanned by  $\{\omega_j\}$.
A \,$2$-from satisfying all the above is called a {\rm (Klein's) fundamental 2-form of second kind}. 
} %%%%%%%%%%%
\footnotetext{We do not use $\partial$ since $x$  and  $y$  are related.}
\vskip 2pt
\fp
{\it Proof.} \ 
Under assuming existence of  $\{\eta_j\}$, 
we see the differential (2.8) satisfies the condition on poles.  
Indeed, we see that, regarding (2.8) as a function of $(x,y)$, 
it has only pole at  $(x,y)=(z,w)$  by 1.15;
and that it is similar as a function of  $(z,w)$  by (2.9).  
A fundamental 2-forms of second kind is obtained similarly 
as \cite{{\bf BG},\,pp.3617--3618} (see also \cite{{\bf Ba},\,around p.194}). 
\qed
\vskip 8pt

It is easily seen that the $\eta_j$ above is written as
  $$
  \aligned
  &\eta_j(x,y)=\frac{h_j(x,y)}{\frac{\partial}{\partial y}f(x,y)}dx \ 
  \ \ \text{\it where}\ h_j(x,y)\in\bold Q[\mu_1,\cdots, \mu_{12}][[x,y]],\ \hbox{\it and} \\ 
  &\hskip 0pt \text{\it $h_j$  is of homogeneous weight.}
  \endaligned
  \tag{2.10}
  $$
Now we finally define  $\eta_j$  uniquely by requiring (see \cite{{\bf BG},\,p.3618})
  $$
  \text{\it the number of terms in  $h_j(x,y)$  is as minimal as possible we could.}  
  \tag{2.11}
  $$
While it is possible to write down the explicit form of  $\eta_j$s, 
we do not need such the expressions in this paper.  
Under the situation above, if we write the 2-form as
  $$
  R\big((x,y),(z,w)\big)
 =\frac{F\big((x,y),(z,w)\big)}
       {(x-z)^2\frac{\partial}{\partial y}f(x,y)\frac{\partial}{\partial w}f(z,w)},
  \tag{2.12}
  $$
we see that  $F\big((x,y),(z,w)\big)$  is a polynomial of homogeneous weight (weight $-8$). 
It is easily seen that  $F\big((x,y),(z,w)\big)$  is 
a polynomial of  $\lambda_j$s, $x$, $y$, $z$, $w$  of homogeneous Sato weight  $-8$.  
We define the period matrices of  $\{\eta_j\}$  by
  $$
  [\,\eta'  \ \ \eta''  ]=\bigg[\int_{\alpha_i}\eta_j \ \ \int_{\beta_i}\eta_j\bigg]_{i,j=1,2,3}. 
  \tag{2.13}
  $$
We concatenate this with (2.3) as
  $$
  M=\left[\matrix\omega' & \omega'' \\ \eta' & \eta'' \endmatrix\right].
  \tag{2.14}
  $$ 
Then, $M$  satisfies
  $$
  M\left[\matrix   & -1_3 \\ 1_3  &  \endmatrix\right]\tp{M}
  =2\pi\sqrt{-1}\left[\matrix   & -1_3 \\ 1_3  &  \endmatrix\right]
  \tag{2.15}
  $$
(see \cite{{\bf Ba}, p.97(c)}, \cite{{\bf FK}, Chap.III}, \cite{{\bf BEL1}, p.11, (1.15); Lemma 2.0.1}). 
This is the {\it generalized Legendre relation} (set of the {\it Weierstrass relations}). 
Especially, we see  ${\omega'}^{-1}\omega''$   is a symmetric matrix.  
It is well-known that 
  $$
   \text{Im}\,({\omega'}^{-1}\omega'') \ \ \ \text{is positive definite}
  \tag{2.16}
  $$
(see \cite{{\bf FK}, Chap.III} for instance). 
%$B=>$D$F(B
%  $$
%  \align
%  \omega'\tp{\omega''}-\omega''\tp{\omega'}&=0, \\ 
%  \eta'\tp{\omega''}-\eta''\tp{\omega'}&=-\tfrac12\pi\sqrt{-1}1_3, \\
%  \eta'\tp{\eta''}-\eta''\tp{\eta'}&=0
%  \endalign
%  $$
%$B$,@.$jN)$D(B. 
It is known by (1.2) that the canonical divisor class is represented by  $4\infty$. 
Hence any theta characteristic is an element of  $\big(\tfrac12\bold Z\big)^{6}$.    
The theta characteristic giving the Riemann constant for our case, namely, taking 
the base point to be  $\infty\in C$  is 
  $$
  \delta=\Big[\matrix\delta'\ \\ \delta''\endmatrix\Big]\in \big(\tfrac12\bold Z\big)^{6}
  \tag{2.17}
  $$
with respect to  $[\,\omega'\ \omega'']$
(see \cite{{\bf Mu}}, pp.163--166, or \cite{{\bf BEL1}}, p.15, (1.18)).  
Under the preparation above, 
we define\footnote{If we redefine $\sigma(u)$  by using another fundamental 2-form of
second kind of 2.5 different from that fixed by (2.11), the  $\sigma(u)$  is 
changed only on the exponential pre-factor. }
  $$
  \aligned
  &\hskip -10pt
  \sigma(u):=\sigma(u;M)=\sigma(u_1,u_2,u_3;M) \\
  &=c\ \text{exp}(-\tfrac{1}{2}u\eta'{\omega'}^{-1}\ ^t\negthinspace u)
  \vartheta\negthinspace
  \left[\delta\right]({\omega'}^{-1}\ ^t\negthinspace u;\ {\omega'}^{-1}\omega'') \\
  &=c\ \text{exp}(-\tfrac{1}{2}u\eta'{\omega'}^{-1}\ ^t\negthinspace u) \\
  &\hskip -10pt\times\sum_{n \in \bold Z^3} \exp \big[2\pi\sqrt{-1}\big\{
  \tfrac12 \ ^t\negthinspace (n+\delta'){\omega'}^{-1}\omega''(n+\delta') 
  + \ ^t\negthinspace (n+\delta')({\omega'}^{-1}u+\delta'')\big\}\big],
  \endaligned
  \tag{2.18}
  $$
where
  $$
  c=\bigg(\frac{\pi^3}{|\omega'|D}\bigg)^{\hskip -2pt\frac12}. 
  \tag{2.19}
  $$
The series in (2.18) converges by (2.16). 
Here  $D$  is the discriminant defined by (1.17) and (1.18), 
$\pi=3.1415\cdots$, 
and  $|\omega'|$  is the determinant of the period matrix  $\omega'$  defined in (2.23). 
The square root of (2.19) is fixed in 2.27 latter. 
In this paper, for  $u\in\bold C^3$, we denote by  $u'$  and  $u''$  
the unique vectors in  $\bold R^3$  such that 
  $$
  u=u'\omega'+u''\omega''.
  \tag{2.20}
  $$
Then we know 
  $$
  \aligned
  L(u,v)&:=\tp{u}(\eta'v'+\eta''v''), \\ 
  \chi(\ell)&:=
  \exp\big\{2\pi\sqrt{-1}\big(\tp{\ell'}\delta''-\tp{\ell''}\delta'+\tfrac12\tp{\ell'}\ell''\big)\big\}
  \ (\in \{1, -1\})
  \endaligned
  \tag{2.21}
  $$
for  $u$, $v\in\bold C^3$  and for  $\ell$ ($=\ell'\omega'+\ell''\omega''$) $\in\Lambda$, 
which is shown by modifying the translational relation of Riemann theta series. 

The following properties are quite important:
\vskip 5pt
\fp
\lbox{84}{360}%--------------------------------------------------------------------------s
{{\bf Lemma 2.22} \ For all  $u\in\bold C^3$,  $\ell\in\Lambda$, 
and  $\gamma\in\text{\rm Sp}(6,{\bold Z})$, we have the following\,{\rm :}\fp
 (1) \ $\sigma(u+\ell;M)=\chi(\ell)\sigma(u;M)\exp L(u+\tfrac12\ell,\ell)$, \fp
 (2) \ $\sigma(u;\gamma M)=\sigma(u;M)$\fp
 (3) \ $u\mapsto\sigma(u;M)$  has zeroes on  $\Theta\br{2}$  of order $1$, \fp
 (4) \ $\sigma(u;M)=0 \iff u\in\Theta\br{2}$. \fp
}%----------------------------------------------------------------------------------------e
\vskip 3pt
\fp
{{\it Proof}. 
The assertion (1) is a special case of  \cite{{\bf Ba}}, p.286, $\ell$.22. 
The assertion (2) is shown by investigating how the transform of  $M$  by  $\gamma$  
corresponds to a change of paths  $\alpha_j$s and $\beta_j$s in (2.2) of period integrals. 
For details, see \cite{{\bf BEL1},\,pp.10--15}.  
The statements (3) and (4) is described in \cite{{\bf Ba}, p.252}}, 
and partially in \cite{{\bf BEL1}}, p.12, Theorem 1.1 and p.15. 
\qed
\vskip 10pt
\fp
{\bf Remark 2.23} \
Let  $M$  is a matrix satisfying (2.15) and (2.16). 
Since the Pfaffian of the Riemann form given by  $L(\ ,\ )$  is $1$ 
as is seen similarly to \cite{{\bf \^O1},\,Lemma 3.1.2}, 
we see that the functions satisfying the equation 2.22(1) form a one dimensional  space 
by \cite{{\bf L}, p.93, Th.4.1}.  
Any such non-trivial solution has properties 2.22(2),\,(3),\,and (4). 
Namely, Lemma 2.22 characterizes the sigma function up to a non-zero multiplicative constant. 
\vskip 10pt
\par
Using 2.22(1),\,(3),\,(4), we have the following equality
  $$
  \aligned
  &\left[
    \frac{\sigma\bigg(\dint_{\infty}^{(x,y)}\omega-\sum_{i=1}^3\dint_{\infty}^{(x_i,y_i)}\omega\bigg)}
         {\sigma\bigg(\dint_{\infty}^{(x,y)}\omega-\sum_{i=1}^3\dint_{\infty}^{(z_i,w_i)}\omega\bigg)}
  \right]
  \left[
  \frac{\sigma\bigg(\dint_{\infty}^{(z,w)}\omega-\sum_{i=1}^3\dint_{\infty}^{(x_i,y_i)}\omega\bigg)}
       {\sigma\bigg(\dint_{\infty}^{(z,w)}\omega-\sum_{i=1}^3\dint_{\infty}^{(z_i,w_i)}\omega\bigg)}
  \right]^{-1}   \\
 &\hskip -20pt =\exp\bigg[
   \int_{(z,w)}^{(x,y)}\bigg(\sum_{i=1}^3\int_{(z_i,w_i)}^{(x_i,y_i)}R\big((x,y),(z,w)\big)dz\bigg)dx\bigg]
   \endaligned
  \tag{2.24}
  $$
by a similar method in \cite{{\bf BEL1}, p.36}. 
We define 
  $$
  \wp_{ij}(u)=-\frac{\partial^2}{\partial u_i\partial u_j}\log \sigma(u), \ \ 
  \wp_{ijk}(u)=-\frac{\partial^3}{\partial u_i\partial u_j\partial u_k}\log \sigma(u), \ \
  \cdots.
  \tag{2.25}
  $$
By Lemma 2.22(1), these are periodic functions with the periods lattice  $\Lambda$. 
Therefore, we can regard these functions as functions on  $J$.  
By (2.24), we see that 
   $$
  \aligned
  &\sum_{i=1}^3\sum_{j=1}^3
   \wp_{ij}\bigg(\dint_{\infty}^{(x,y)}\omega-\sum_{r=1}^3\dint_{\infty}^{(x_r,y_r)}\omega\bigg)
   \frac{\omega_i(x,y)}{dx}\frac{\omega_j(x_r,y_r)}{dx_r} \\
 &=\frac{F\big((x,y),(x_r,y_r)\big)}
        {(x-x_r)^2\frac{\partial}{\partial y}f(x,y)\frac{\partial}{\partial y_r}f(x_r,y_r)}.
   \endaligned
   \tag{2.26}
   $$
Taking a local parameter  $t$  at  $\infty$  such that  $x=1/t^3$, 
and expanding both sides of (2.26) as power series of  $t$  as in \cite{{\bf BG}, p.3621}, 
we know that  $\wp_{j_1j_2\cdots j_n}(u)$s are of homogeneous weight.  
Therefore,  $\sigma(u)$  is of homogeneous weight
\footnote{Consider the product of 2 functions.  
If one of them is not homogeneous, then the product of them is also 
not homogeneous. We use this fact here.}.  
Including this fact, we see the following. 
\vskip 7pt
\fp
\lbox{91}{360}%------------------------------------------------------------------------------s
{{\bf Lemma-Definition 2.27} \ \it The power series expansion of  $\sigma(u)$  at  $u=(0,0,0)$  
with respect to  $u_1$, $u_2$, $u_3$  has of homogeneous Sato weight  $5$, 
and it is of the form 
  $$
  \sigma(u)=\pm\big(u_1-u_3{u_2}^2+\tfrac1{20}{u_3}^5\big)
           +(d^{\circ}({\lambda_1,\cdots,\lambda_{12}})\geqq 1).
  $$
We define precisely  $\sigma(u)$  by taking  $+$  in the  $\pm$  above.  
}%----------------------------------------------------------------------------------------e
\vskip 3pt
\fp
{\it Proof}. We see that  $\sigma(u)\in\bold Q[\lambda_1,\cdots, \lambda_{12}][[u_1,u_2,u_3]]$  
by as described above. 
The homogeneousness in Sato weight also is already explained. 
We know the first three terms up to a multiplicative non-zero constant 
from the main result and Example 5.5 of \cite{{\bf BEL2}}. 
While such the constant is known by the method of Frobenius (\cite{{\bf Fr}})
by using a generalization of Thomae's formula, 
we omit the proof because of less importance in the below.  
\qed
\vskip 10pt
\fp
\lbox{53}{360}{{\bf Lemma 2.28} \  \it The function  $\sigma(u)$  is an odd function. 
Namely, we have
  $$
  \sigma(\lceil -1\rceil u)=-\sigma(u).
  $$
}
\vskip 3pt
\fp
{\it Proof}. \ 
For the period lattice, we have  $\lceil -1\rceil \Lambda=\Lambda$. 
This is seen by considering the integral of  $\ell\in\Lambda$  in the opposite directions. 
Adding this with 2.23, we see that there is a constant  $K$  such that 
  $$
  \sigma(\lceil -1\rceil u)=K\sigma(u).
  \tag{2.29}
  $$
%  $$
%  [\zeta]u\du{1}=\zeta^{\nu} u\du{1}
%  \tag{5.6}
%  $$
Since the weight 5 is odd integer, Lemma 2.27 implies  $K=-1$. 
\qed
\vskip 10pt
%\fp
%{\gt\bf $BCm0U(B 2.28}  \ $\sigma(u)$  $B$Nf25i?tE83+$O(B, $B62$i$/(B
%$\bold Z[\lambda_1,\cdots,\lambda_s]$  $B>e$N(B  $u\du{w_g}$, $\cdots$, $u\du{w_1}$  
%$B$K4X$9$k(B Hurwitz $B@0$Jf25i?t$G$"$i$&$H;W$O$l$k(B. 
%\fp
%(2)  $\{\wp_{ij}(u)\,|\, i\geqq j,\ 1\leqq i\leqq 3,\ 1\leqq j\leqq 3\}$  $B$,@~7?FHN)$G$"$k$3$H$O(B, 
%$BJdBj(B 2.27 $B$NE83+$N:G=i$N$$$/$D$+9`$+$i3N$+$a$k$3$H$,$G$-$k(B. 
%Jacobi $BB?MMBN(B  $J$  $B$O(B  $3$  $B<!85$J$N$G(B, $B$3$l$i(B  $3(3+1)/2=6$  $B8D$NH!?t$O(B  $J$  $B$N(B
%$BH!?tBN$r@8@.$9$k$3$H$,$o$+$k(B.  
\par
In the sequell, we write simply 
  $$
  \sigma_j(u)=\frac{\partial}{\partial u_j}\sigma(u), \ \ \ 
  \sigma_{ij}(u)=\frac{\partial^2}{\partial u_i\partial u_j}\sigma(u).
  \tag{2.30}
  $$
\vskip 10pt
\fp
\lbox{133}{360}%------------  ($BJdBj(B) -------------------------------------------------------s
{{\bf Lemma 2.31} \ \it Let  $u$, $u\lr{1}$, $u\lr{2}$, $v\in\kappa^{-1}\iota(C)$. 
Then we have the following\,{\rm :} \fp
{\rm (1)}  $\sigma(u\lr{1}+u\lr{2})=0$\,{\rm ;} \fp
{\rm (2)}  The expansion of  $v\mapsto\sigma(u\lr{1}+u\lr{2}+v)$  with respect to  $v_3$  
is of the form\fp
\hskip 23pt $\sigma(u\lr{1}+u\lr{2}+v)=\sigma_3(u\lr{1}+u\lr{2})v_3+(d^{\circ}(v_3)\geqq 2)$\,{\rm ;}\fp
{\rm (3)}  $\sigma_3(u)=0$\,{\rm ;} \fp
{\rm (4)}  The expansion of  $v\mapsto\sigma_3(u+v)$  with respect to  $u_3$  is of the form\fp
\hskip 23pt $\sigma_3(u+v)=\sigma_{33}(u)v_3+(d^{\circ}(v_3)\geqq 2)$\,{\rm ;} \fp
{\rm (5)}  The expansion of  $v\mapsto\sigma_{33}(v)$  with respect to  $v_3$  is of the form\fp
\hskip 23pt $\sigma_{33}(v)={v_3}^3+(d^{\circ}(v_3)\geqq 4)$.  
}
%----------------------------------------------------------------------------------------e
\vskip 3pt
\fp
{\it Proof}. The assertions (1) and (2) are repetition of 2.22(3).  
For the expansion of (2), 
%  $$
%  \sigma(u\lr{1}+u\lr{2}+v)=\sigma_3(u\lr{1}+u\lr{2})v_3+(d^{\circ}(v_3)\geqq 2), 
%  \tag{2.32}
%  $$
by taking   $u\lr{2}$  close to  $(0,0,0)$, we see
\footnote{In this situation, for given point  $P\in C$ ($\neq \infty$),
the rank of the Brill-Noether matrix  $B(P+\infty)$  
of the divisor  $P+\infty$  on  $C$  is  $\rank B(P+\infty)=2$, and  $\dim\Gamma(C,\Cal O(P))=1$. 
It is impossible to show directly  $\sigma_3(u)=0$ ($u\in\kappa^{-1}\iota(C)$) 
from Riemann singularity theorem. }
$\sigma_3(u\lr{1})=0$  
because of  $\sigma(u\lr{1}+v)=0$  by (1). 
The assertion (4) is obviously follows from (3).  
The assertion (5) is seen from 2.27.  \qed
\vskip 10pt
\fp
\lbox{54}{360}%-------------------------------------------------------------------------s
{{\bf Lemma 2.32} \ \it We have the following translational relations\,{\rm :}
\fp
{\rm (1)}  For  $u\in\kappa^{-1}(\Theta\br{2})$, we have 
     $\sigma_3(u+\ell)=\chi(\ell)\sigma_3(u)\exp L(u+\tfrac12\ell,\ell)$\,{\rm ;} \fp
{\rm (2)}  For  $u\in\kappa^{-1}(\Theta\br{1})$, we have  
     $\sigma_{33}(u+\ell)=\chi(\ell)\sigma_{33}(u)\exp L(u+\tfrac12\ell,\ell)$. 
}%---------------------------------------------------------------------------------------e
\vskip 3pt
\fp
{{\it Proof}. Differentiating both side of 2.22(1) by  $u_3$  once or twice, 
we see the assertions by using 2.31. \qed
\vskip 10pt
\fp
\lbox{43}{360}%-------------------------------------------------------------------------s
{{\bf Lemma 2.33} \ \it The function  $u\mapsto\sigma_{33}(u)$  on  $\kappa^{-1}\iota(C)$  
has only zero of order  $3$  at  $u=(0,0,0)$  modulo  $\Lambda$,
and no zeroes elsewhere. }
\vskip 3pt
\fp
{{\it Proof}. After taking logarithm of the equation in 2.33(2),
by integrating it around a regular polygon of the Riemann surface given by  $C$,  
the principle of arguments and the general Legendre relation, as usual, 
shows the total of orders of zeroes of the claimed function is $3$.   
Hence, 2.31(5) implies the assertion. \qed

\newpage
%@@@@@@@@@@@@@@@@@@@@@@@@@@@@@@@@@@@@@@@@@@@@@@@@@@@@@@@@@@@@@@@@@@@@@@@@@@@@@@@@@@@@@@@@@@@
\fp
{\twelveptbf 3. The sigma functions (Purely trigonal case)} 
\vskip 5pt
\fp
The rest of this paper, we treat only the curve  $C$ 
  $$
  y^3 =x^4+\lambda_3x^3+\lambda_6x^2+\lambda_9x+\lambda_{12} \ \ \text{($\lambda_j$  are constants)}
  \tag{3.1}
  $$
that specializing (1.1).  Then (1.2) is given by
  $$
  \omega_1=\frac{dx}{3y^2}, \ \ 
  \omega_2=\frac{xdx}{3y^2}, \ \ 
  \omega_2=\frac{ydx}{3y^2}=\frac{dx}{3y};
  \tag{1.$2'$}
  $$
and (2.4) is 
  $$
  \Omega\big((x,y),(x',y')\big)=\frac{y^2+yy'+{y'}^2}{(x-x')3y^2}.
  \tag{2.$4'$}
  $$
Using these, we define  $\sigma(u)$  as in the Section 2. 
Let  $\zeta=e^{2\pi \sqrt{-1}/3}$. 
Since the curve  $C$  has an automorphism defined by  $(x,y)\mapsto (x, {\zeta}y)$, 
${\zeta}^j$  acts for  $u=(u_1,u_2,u_3)\in\kappa^{-1}\iota(C)$  by
  $$
  {\lceil{\zeta}^j\rceil}u=({\zeta}^ju_1,{\zeta}^ju_2,{\zeta}^{2j}u_3)
                        =\int_{\infty}^{(x,{\zeta}^jy)}\omega.
  \tag{3.2}
  $$
Hence, we see that  $W\br{k}$  and  $\Theta\br{k}$  are stable under the action by  ${\zeta}^j$. 
Under this situation, the function  $\sigma(u;M)$  has the following special property. 
\vskip 5pt
\fp
\lbox{44}{360}%----------------------------------------------------------------------------s
{{\bf Lemma 3.3} \ \it We have
  $$
  \sigma({\lceil\zeta\rceil}u)=\zeta\sigma(u). 
  $$
}%-----------------------------------------------------------------------------------------e
\vskip 0pt
\fp
{\it Proof}. \  
This follows from 2.23, the fact  $\lceil\zeta\rceil\Lambda=\Lambda$, and 2.27.  
\qed
\vskip 10pt

The following Lemma is the key for our main result. 
\vskip 3pt
\fp
\lbox{130}{360}
{{\bf Lemma 3.4} \ \it
For  $u\in\kappa^{-1}\iota(C)$, the function
  $$
  v\mapsto \sigma_3(u+{\lceil\zeta\rceil}v), \ \ 
  \big(\text{resp.}\ v\mapsto \sigma_3(u+{\lceil\zeta^2\rceil}v)\big)
  $$
on  $\kappa^{-1}\iota(C)$  has 
only $3$  zeroes  $(0,0,0)$, $u$, ${\lceil\zeta^2\rceil}u$  
{\rm (}resp. ${\lceil\zeta\rceil}u$ {\rm resp.)} of order  $1$  
modulo  $\Lambda$, and has no zeroes elsewhere. 
Its power series expansion at  $(0,0,0)$  with respect to  $u_3$  is of the form
  $$
  \align
  \sigma_3(u+{\lceil\zeta\rceil}v)&=\sigma_{33}(u)v_3+(d^{\circ}(v_3)\geqq 2)\\
  \big(\text{resp.}  \ 
  \sigma_3(u+{\lceil\zeta^2\rceil}v)&=\sigma_{33}(u)v_3+(d^{\circ}(v_3)\geqq 2)\big).
  \endalign
  $$
}% ----------------------------------------------------------------------------------------e
\vskip 0pt
\fp
{\it Proof}. \ 
Remark 2.28 shows that  $u\mapsto\sigma_3(u)$  is an even function. 
Becase  $u+{\lceil\zeta\rceil}u+{\lceil\zeta^2\rceil}u=(0,0,0)$, 
we have for $u\in\kappa^{-1}\iota(C)$ 
  $$
  \sigma_3(u+{\lceil\zeta\rceil}u)=\sigma_3(-{\lceil\zeta^2\rceil}u)
  =\sigma_3({\lceil\zeta^2\rceil}u)=\zeta^2\sigma_3(u)=0
  \tag{3.5}
  $$
by 2.31(3) and 3.3.  
The remained assertions follows from 2.31(4) and 3.3. \qed
\vskip 10pt
\par
The following Lemma is used in Lemma 5.1. 
\vskip 3pt
\fp
\lbox{59}{360}
{{\bf Lemma 3.6} \ \it For $u\in\kappa^{-1}\iota(C)$, we have\,{\rm :} 
  $$
  \frac{\sigma_3(2u)}{\sigma_{33}(u)^4}=3y(u)^2. 
  $$
}%----------------------------------------------------------------------------------------e
\vskip 5pt
\fp
{\it Proof}. \ Because of 2.27 and 1.10, we have that 
  $$
  \sigma_3(2u)=\tfrac14(2u_3)^4-\big(2(\tfrac12{u_3}^2+\cdots)\big)^2+\cdots 
              =3{u_3}^4+\cdots,
  $$
so that 
  $$
  \frac{\sigma_3(2u)}{\sigma_{33}(u)^4}=\frac{3{u_3}^4+\cdots}{\big({u_3}^3+\cdots\big)^4}
                                       =\frac3{{u_3}^8}+\cdots. 
  \tag{3.7}
  $$
The left hand side of this is a periodic function 
with respect to  $\Lambda$  by 2.32.  
Therefore, its pole is only at  $(0,0,0)$  modulo  $\Lambda$  by 2.31. 
By 3.3 we see 
  $$
  \frac{\sigma_3({\lceil\zeta\rceil}2u)}{\sigma_{33}({\lceil\zeta\rceil}u)^4}
  =\zeta^2\frac{\sigma_3(2u)}{\sigma_{33}(u)^4}. 
  $$
and that the right hand side of (3.7) should be  $3y(u)^2$.  
\qed

\newpage

%%%%%%%%%%%%%%%%%%%%%%%%%%%%%%%%%%%%%%%%%%%%%%%%%%%%%%%%%%%%%%%%%%%%%%%%
\fp
{\twelveptbf 4. Frobenius-Stickelberger-Type Fromulae} 
\vskip 3pt
\fp
The initial case of Frobenius-Stickelberger type formula for 
the curve  $C$\,:\,$y^3=x^4+\cdots$  is as follows:
\vskip 5pt
\fp
\lbox{63}{360}%---------------------------------------------------------------------------s
{{\bf Proposition 4.1} \it \ For  $u$, $v\in\kappa^{-1}\iota(C)$  we have\,{\rm :}
  $$
  \frac{
    \sigma_3(u+v)
   \,\sigma_3(u+{\lceil\zeta\rceil}v)
   \,\sigma_3(u+{\lceil{\zeta}^2\rceil}v)}
   {\sigma_{33}(u)^3
    \sigma_{33}(v)^3} 
   =\big((x(u)-x(v)\big)^2 
   ={\eightpoint 
  \left|\matrix
   1 &  x(u) \\
   1 &  x(v) \\
   \endmatrix\right|^2}.
  $$
}%----------------------------------------------------------------------------------------e
\vskip 5pt
\fp
{\it Proof}. \ 
Lemma 2.33 shows the left hand side is a periodic function of  $v$ (resp. $u$) with the periods  $\Lambda$.  
Now,we regards the left hand side as a function of  $v$.    
Lemma 3.4 states that the second and third factors vanish at  $v=u$  modulo  $\Lambda$,  
namely the left hand side has zero of order 2 there, 
and has two zeroes of order 1 at  ${\lceil\zeta\rceil}u$  and  ${\lceil\zeta^2\rceil}u$  modulo  $\Lambda$.  
It has no zero elsewhere.  
Its only pole is  $v=(0,0,0)$  modulo  $\Lambda$  by Lemma 2.34.  
The order is  $3\times 3-3=6$.  
These situations are exactly the same for the right hand side.  
Therefore the two sides coincide up to a multiplicative constant depending only on  $u$.  
If we expand both sides with respect to  $v_3$  
we see that the coefficients of the least terms of the two sides coincide exactly
by using Lemmas 2.31 and 1.15. Thus, the proof has completed.  
\qed
\vskip 10pt
\fp
{\bf Remarks 4.2} \ 
Note that  $-u\not\in\kappa^{-1}\iota(C)$  in general.  
Why the initial formula above different from the initial formula for 
the case of hyperelliptic functions is explained by this fact.  \fp
\vskip 10pt
\fp
\lbox{194}{390}%===============================================================================s
{{\bf Theorem 4.3}\,{\rm (Frobenius-Stickelberger type formula)} \it  
Let  $n\geqq 3$  be an integer.  
Let  $\sigma(u)$, $x(u)$, and  $y(u)$  are those defined 
for the curve  $C$ {\rm :} $y^3=x^4+\cdots$  as above. 
Assume  $u\lr{1}$, $\cdots$, $u\lr{n}$  are points on  $\kappa^{-1}\iota(C)$.  
Then we have\,{\rm :}
  $$
  \align
  &\frac{
  \sigma(u\lr{1}+u\lr{2}+\cdots+u\lr{n})
  \prod_{i<j}\sigma_3(u\lr{i}+{\lceil\zeta  \rceil} u\lr{j})\,
                  \sigma_3(u\lr{i}+{\lceil\zeta^2\rceil} u\lr{j})
  }{
  \sigma_{33}(u\lr{1})^{2n-1}\cdots\sigma_{33}(u\lr{n})^{2n-1}
  }\\
  &=
  {\eightpoint
  \left|
  \matrix
  1  &  x(u\lr{1})  &  y(u\lr{1})    &  x^2(u\lr{1})  
                    &  yx(u\lr{1})   &  y^2(u\lr{1})  &  x^3(u\lr{1})  
                    &  yx^2(u\lr{1}) &  y^2x(u\lr{1}) &  \cdots   \\
  1  &  x(u\lr{2})  &  y(u\lr{2})    &  x^2(u\lr{2})  
                    &  yx(u\lr{2})   &  y^2(u\lr{2})  &  x^3(u\lr{2})  
                    &  yx^2(u\lr{2}) &  y^2x(u\lr{2}) &  \cdots   \\
  \vdots &  \vdots  &  \vdots        &  \vdots  
                    &  \vdots        &  \vdots        &  \vdots  
                    &  \vdots        &  \vdots        &  \ddots  \\
  1  &  x(u\lr{n})  &  y(u\lr{n})    &  x^2(u\lr{n})  
                    &  yx(u\lr{n})   &  y^2(u\lr{n})  &  x^3(u\lr{n})  
                    &  yx^2(u\lr{n}) &  y^2x(u\lr{n}) &  \cdots   \\
  \endmatrix\right|
  } \\
  &\hskip 20pt
  {\eightpoint
  \cdot
  \left|
  \matrix
  1  &  x(u\lr{1})  &  x^2(u\lr{1})    &  \cdots      &  x^{n-1}(u\lr{1})   \\
  1  &  x(u\lr{2})  &  x^2(u\lr{2})    &  \cdots      &  x^{n-1}(u\lr{2})   \\
  \vdots &  \vdots  &  \vdots          &  \ddots      &  \vdots             \\
  1  &  x(u\lr{n})  &  x^2(u\lr{n})    &  \cdots      &  x^{n-1}(u\lr{n})   \\
  \endmatrix
  \right|.}
  \endalign
  $$
}%=====================================================================================e
\vskip 10pt
\fp
{\bf Remark 4.4} \ The each row of the first determinant consists of monomials of  $x$  and  $y$ 
(they have a pole only at  $(0,0,0)$)  displayed according to their order of the pole,  
and their order are 
  $$
  0, \ 3, \ 4, \ 6, \ 7, \ 8, \ 9, \ \cdots
  $$
respectively.
%\fp
%(2) We see that if the second determinant also had contained a few column of  
%monomial including  $y(u)$, then the right hand side could never vanish 
%at  $v={\lceil\zeta\rceil}u$, but in fact  $v\mapsto \sigma_{33}(u+{\lceil\zeta\rceil}v)$  
%has a zeroe at  ${\lceil\zeta\rceil}u$.
\vskip 10pt
\fp
{\it Proof}. \ We prove the formula by induction.  
First of all, we see by Lemma 2.33 that the left hand side is a periodic function of  $u$  and of  $v$  
with the periods  $\Lambda$. \fp  
(1) For the case  $n=3$, the reader easily prove the formula, 
if he is in mind the proof of the case  $n\geqq 4$  below (and \cite{{\bf \^O3}}, p.309), 
and so we omit the proof of this case.  \fp
\fp
(2) Suppose  $n\geqq 4$. 
We regards both sides as functions of  $u:=u\lr{n}$.  
\fp
(2-a) We know the divisor of the two sides by Lemmas 2.34, 2.22(3),(4), 3.4 as follows:
\fp
{\it The left hand side}. 
The numerator of the left hand side has zeroes of order 2 at  $(n-1)$  points
  $$
  u\lr{j}\ \text{modulo  $\Lambda$} \ \ \ 
  (j=1,\,\cdots,\,n-1),
  $$
zeroes of order 1 at  $2(n-1)$  points
  $$
  {\lceil\zeta  \rceil}u\lr{j}, \ \ 
  {\lceil\zeta^2\rceil}u\lr{j}\ \text{modulo  $\Lambda$} \ \ \ 
  (j=1,\,\cdots,\,n-1),
  $$
and zeroes of order $2(n-1)$  at  $u=(0,0,0)$.  
It has no zeroes elsewhere.  
The denominator has only zero at  $u=(0,0,0)$  of order  $3(2n-1)=6n-3$.  
Therefore the left hand side has only pole at  $(0,0,0)$  
of order  $(6n-3)-2(n-1)=4n-1$.  
There are  $(4n-1)-2(n-1)-2(n-1)=3$  unknown zeroes of the left hand side. 
\fp
{\it The right hand side}.
The situation is exactly the same in this side.  
Indeed, both determinants vanish at 
  $$
  u\lr{j}\ \text{modulo  $\Lambda$} \ \ \ (j=1,\,\cdots,\,n-1),
  $$
and only the second determinant vanishes at 
  $$
  {\lceil\zeta  \rceil}u\lr{j}, \ \
  {\lceil\zeta^2\rceil}u\lr{j}\ \text{modulo  $\Lambda$} \ \  \ 
  (j=1,\,\cdots,\,n-1).  
  $$
The deepest pole comes from  $(n,n)$-entries of the two determinants.  
The sum of the pole orders is  $(n-2)+3(n-1)=4n-1$.  
\fp
(2-b) We let  $\alpha$, $\beta$, $\gamma\in\bold C^3$ modulo $\Lambda$  
are the unknown zeroes of the right hand side. 
Then the Abel-Jacobi theorem shows
  $$
  \align
  2(u\lr{1}+\cdots+u\lr{n-1})
  +({\lceil\zeta  \rceil}u\lr{1}\cdots+{\lceil\zeta  \rceil}u\lr{n-1})&
  +({\lceil\zeta^2\rceil}u\lr{1}\cdots+{\lceil\zeta^2\rceil}u\lr{n-1}) \\
  &+\alpha+\beta+\gamma\in\Lambda.
  \endalign
  $$
Since  $u\lr{j}+{\lceil\zeta\rceil}u\lr{j}+{\lceil\zeta^2\rceil}u\lr{j}=(0,0,0)$, 
we see
  $$
  u\lr{1}+\cdots+u\lr{n-1}
  +\alpha+\beta+\gamma\in\Lambda, 
  $$
so that 
  $$
  u\lr{1}+\cdots+u\lr{n-1}+u\lr{n}
  \equiv u-\alpha-\beta-\gamma \mod{\Lambda}.
  $$
Observing the first factor in the numerator of the left hand side, 
we see the left hand side has zeroes of order 1 at  $\alpha$, $\beta$, $\gamma$, too. 
Thus the two sides coincide up to multiplicative constant.  
\fp
(2-c) The coefficients of the lowest term with respect to Laurent expansion with respect to  $u$,
is just the hypothesis of induction, and they coincide exactly.  
Hence the proof has completed. 
\qed

\newpage

\fp
{\twelveptbf 4. Kiepert Type Formula}
\vskip 5pt
\fp
We prove the following fundamental formula.
\vskip 3pt
\fp
\lbox{58}{360}%-----------------------------------------------------------------------------s
{{\bf Lemma 5.1.} \ \it We have 
  $$
  \lim_{v\to u}
  \frac{\sigma_3(u+{\lceil\zeta\rceil}v)\,\sigma_3(u+{\lceil\zeta^2\rceil}v)}
       {\sigma_{33}(u)\sigma_{33}(v)(u_3-v_3)^2}
  =3x(u)^2.
  $$
}%-------------------------------------------------------------------------------------------e
\vskip 5pt
\fp
{\it Proof}. \ Lemma 3.6 and Proposition 4.1 show
  $$
  \align
  3y(u)^2
  \bigg(\lim_{v\to u}&
  \frac{\sigma_3(u+{\lceil\zeta\rceil}v)\,\sigma_3(u+{\lceil\zeta^2\rceil}v)}
       {\sigma_{33}(u)\sigma_{33}(v)(u_3-v_3)^2}
  \bigg) \\
  &=\lim_{v\to u}
  \frac{
    \sigma_3(u+v)
  \,\sigma_3(u+{\lceil\zeta\rceil}v)
  \,\sigma_3(u+{\lceil{\zeta}^2\rceil}v)}
   {\sigma_{33}(u)^3
    \sigma_{33}(v)^3
    (u_3-v_3)^2} \\
  &=\lim_{v\to u}\Big(\frac{x(u)-x(v)}{u_3-v_3}\Big)^2 \\
  &=\lim_{v\to u}\Big(\frac{dx}{du_3}(u)\Big)^2 \\
  &=\Big(\frac{3y(u)}{x(u)}\Big)^2.
  \endalign
  $$
This yields our desired formula. 
\qed
\vskip 10pt
\fp
\lbox{145}{395}
{{\bf Corollary 5.2.}\ (Kiepert-type formula) \it
Suppose $n\geqq 3$  and  $u\in\kappa^{-1}\iota(C)$. We have 
  $$
  \align
  &\psi_n(u):=\frac{\sigma(nu)}
  {\sigma_{33}(u)^{n^2}}
  =y^{n(n-1)/2}(u)\times \\
  &{\eightpoint
  \left|
  \matrix
  x'         &  y'             &  (x^2)'  
             &  (yx)'          &  (y^2)'         
             &  (x^3)'         &  (yx^2)'        &  (y^2x)'      
             &  \cdots   \\
  x''        &  y''            &  (x^2)''  
             &  (yx)''         &  (y^2)''        
             &  (x^3)''        &  (yx^2)''       &  (y^2x)''
             &  \cdots   \\
 \vdots      &  \vdots         &  \vdots  
             &  \vdots         &  \vdots         
             &  \vdots         &  \vdots         &  \vdots
             &  \ddots   \\
  x\lr{n-1}  &  y\lr{n-1}      &  (x^2)\lr{n-1}  
             &  (yx)\lr{n-1}   &  (y^2)\lr{n-1}  
             &  (x^3)\lr{n-1}  &  (yx^2)\lr{n-1} &  (y^2x)\lr{n-1} 
             &  \cdots   \\
  \endmatrix\right|(u)},
  \endalign
  $$
where  $'$  means   $\tfrac{d}{du_3}$  and the determinant is of size  $(n-1)\times(n-1)$. 
}
\vskip 5pt
\fp
{\it Proof}. \ 
Each factor  $x(u\lr{i})-x(u\lr{j})$  of the Vandermonde determinant gives 
  $$
  \lim_{v\to u}
  \frac{x(u)-x(v)}{u_3-v_3}=\frac{dx}{du_3}(u)=3\frac{y(u)}{x(u)}.  
  $$
Then Theorem 4.3 and Lemma 5.1 gives the formula of Kiepert type 
by similar manipulation of the proof of \cite{{\bf \^O2}, Theorem 3.3}. \qed

\newpage

%\hyphenation{Entwickelungs-co-{\"e}f-fi-cienten

\enddocument

\Refs\nofrills{\it References}
% \baselineskip=2\baselineskip
\widestnumber\key{[{\bf BEL3}]}
% \widestnumber\no{{\bf 18}}
\ref
  \key   {\bf BEL1}
  \by    {V.M. Buchstaber, V.Z. Enolskii and D.V. Leykin}  
  \paper {\rm Kleinian functions, hyperelliptic Jacobians and applications}
  \jour  {\it Reviews in Math. and Math. Physics }
  \vol   {\bf 10}
  \yr    1997
  \pages 1--125
\endref
\ref
  \key   {\bf BEL2} % \key  BEL
  \by    {V.M. Buchstaber, V.Z. Enolskii, and D.V. Leykin}
  \paper {\rm Rational analogs of Abelian functions}
  \jour  {\it Functional Anal. Appl. }
  \vol   {\bf 33}
  \yr    1999
  \pages 83--94
\endref
\ref
  \key   {\bf BLE3} % \key  BLE3
  \by    {V.M. Buchstaber, D.V. Leykin, and V.Z. Enolskii}
  \paper {\rm $\sigma$-function of $(n,s)$-curves}
  \jour  {\it Russ. Math. Surv.}
  \vol   {\bf 54}
  \yr    1999
  \pages 628-629
\endref
\ref   
  \key   {\bf FS} % \key  FS
  \by    {F.G. Frobenius and L. Stickelberger}
  \paper {\rm Zur Theorie der elliptischen Functionen}
  \jour  {\it J. reine angew. Math.}
  \vol   {\bf 83}
  \yr    1877
  \pages 175--179
\endref
\ref
  \key   {\bf Fr} 
  \by    {Frobenius, G.}
  \paper {\rm \"Uber die constanten Factoren der Thetareihen}
  \pages 244-263
  \yr    1885
  \vol   98
  \jour  {\it J. reine angew. Math. }
\endref
\ref
  \key   {\bf K} % \key  K
  \by    {L. Kiepert}
  \paper {\rm Wirkliche Ausf\"uhrung der ganzzahligen Multiplikation 
         der elliptichen Funktionen} 
  \jour  {\it J. reine angew. Math.}
  \vol   {\bf 76}
  \yr    1873 
  \pages 21--33
\endref
\ref
  \key   {\bf \^O1} % \key  \^O1
  \by    {\^Onishi, Y.}
  \pages 381-431
  \paper {\rm  Complex multiplication formulae for hyperelliptic curves of genus three}
  \yr    1998
  \vol   21
  \jour  {\it Tokyo J. Math. }
\endref
\ref
  \key   {\bf \^O2} % \key  \^O2
  \by    {Y. \^Onishi}
  \pages 353-364
  \paper {\rm Determinant expressions for Abelian functions in genus two}
  \yr    2002
  \vol   44
  \jour  {\it Glasgow Math. J. }
\endref
\ref
  \key   {\bf \^O3} % \key  \^O4
  \by    {Y. \^Onishi}
  \pages 299-312
  \paper {\rm Determinant expressions for hyperelliptic functions in genus three}
  \yr    2004
  \vol   27
  \jour  {\it Tokyo J. Math. }
\endref
\ref
  \key   {\bf \^O4} % \key  \^O4
  \by    {Y. \^Onishi}
  \toappear
%  \pages 299-312
  \paper {\rm Determinant expressions for hyperelliptic functions (with an appendix by Shigeki Matsutani)}
%  \yr    2005
%  \vol   27
  \jour  {\it Proc. Edinburgh Math. Soc.}
\endref
\endRefs
\bye